\documentclass[a4paper,12pt]{article}
\usepackage{amssymb,amsmath}
\newcommand{\ssl}{\mathsf{\Lambda}}
\renewcommand{\b}[1]{{\bf #1}}
\newcommand{\Proj}{\mathbb{P}}
\newcommand{\x}{{\bf x}}
\newcommand{\hcf}{{\rm h.c.f.}}
\newcommand{\ep}{\varepsilon}

\newcommand{\Q}{\mathbb{Q}}
\newcommand{\R}{\mathbb{R}}
\numberwithin{equation}{section}
\newcommand{\N}{\mathbb{N}}
\newcommand{\Z}{\mathbb{Z}}
\newcommand{\cl}[1]{{\cal #1}}
\renewcommand{\mod}[1]{\hspace{-2.9mm}\pmod{#1}}

\newtheorem{theorem}{Theorem}
\newtheorem{lemma}{Lemma}

\begin{document}

\title{The Density of Rational Points on Cayley's Cubic Surface}
\author{D.R. Heath-Brown\\Mathematical Institute, Oxford}
\date{}
\maketitle

\begin{center}
{\bf Abstract}
\end{center}

{\small
The Cayley cubic surface is given by the equation $\sum_{i=1}^4
X_i^{-1}=0$.  We show that the number of non-trivial primitive integer
points of size at most $B$ is of exact order $B(\log B)^6$, as
predicted by Manin's conjecture.}

\section{Introduction}

The Cayley cubic surface is defined in $\Proj^3$ by the equation
\[\frac{1}{X_1}+\frac{1}{X_2}+\frac{1}{X_3}+\frac{1}{X_4}=0\]
or equivalently by
\[C:\; C(X_1,X_2,X_3,X_4)=X_2X_3X_4+X_1X_3X_4+X_1X_2X_4+X_1X_2X_3=0.\]
It has four singularities, at the points $(1,0,0,0),
(0,1,0,0), (0,0,1,0)$ and \linebreak $(0,0,0,1)$.  Moreover there are exactly 9
lines in the surface, and all of these are defined over the
rationals.  Three of the lines have the form $X_i+X_j=X_k+X_l=0$, while 
the remaining six have the shape $X_i=X_j=0$.  We shall write $U$ for
the complement of these lines in the surface $C$.

The aim of this paper is to consider the density of rational points on the
surface $C$.  It transpires that ``most'' of the rational points lie
on one of the lines described above.  We shall think of such points as
being ``trivial'' and exclude them from our counting function.  We
therefore define
\[N^*(B)=\#\{\x\in\Z^4:\;\x\in U, \;\max|x_i|\le B\},\]
where $\x=(x_1,x_2,x_3,x_4)$.  Indeed, since vectors $\x$ which are
scalar multiples of each other represent the same projective point, it
is natural to consider only primitive vectors $\x$.  (A vector $\x$ is
said to be primitive if ${\rm h.c.f.}(x_1,x_2,x_3,x_4)=1$.)  With this
in mind we set
\[N(B)=\#\{\x\in\Z^4: \;\x\in U,\;\max|x_i|\le B,\;\x\;\mbox{primitive}\}.\]
The corresponding number of rational points in $\Proj^3$ is 
$\frac{1}{2}N(B)$, since
$\x$ and $-\x$ represent the same point.  Our two counting functions
are closely related, since 
\[N^*(B)=\sum_{h\le B}N(B/h),\]
as one readily verifies.

Manin (see Batyrev and Manin \cite{Manin}) 
has given a very general conjecture which would
predict in our case that
\[N(B)\sim cB(\log B)^{6},\]
for a suitable positive constant $c$.
For an arbitrary cubic surface one expects something of this type, with
the exponent of the logarithm 
being one less that the rank of the Picard group of
the surface, and in our case this rank is 7.  Unfortunately the
conjecture has only been established for a small number of extremely
simple cubic
surfaces, all of which are singular.  For example, several authors
have considered the surface $X_1X_2X_3=X_4^3$, see de la Bret\`{e}che 
\cite{dlB}, Fouvry \cite{Fouv}, Heath-Brown and Moroz \cite{HBM} and
Salberger \cite{Salb}.  The Cayley surface, while still singular, is
considerably more intricate than any previous example.  The goal of
the present paper is to establish the following estimates.

\begin{theorem}
We have
\[B(\log B)^{6}\ll N(B)\ll B(\log B)^{6}.\]
\end{theorem}

Of the two inequalities here, the lower bound is relatively easy to
prove.  Indeed Slater and Swinnerton-Dyer \cite{SSD} have
established the lower bound corresponding to Manin's conjecture for
any non-singular cubic surface defined over $\Q$, providing that it contains
two skew lines defined over $\Q$.  Although our surface is singular, it
does contain several pairs of skew lines, and these are crucial to our
argument.   It would have been somewhat easier to have established
upper bounds of order $B^{1+\ep}$, with an arbitrary positive constant
$\ep$, or indeed of order $B(\log B)^A$ for some large constant $A$.
However to achieve the correct exponent 6 requires more work.

It is natural to ask how close we come to establishing an asymptotic
formula for $N(B)$.  An analysis of the argument in \S 6 shows that
the difficulty arises through our use of Lemma 6, which gives an upper
bound for the number of primitive lattice points in $\Z^3$, lying in a
box, and which
satisfy a given linear equation.  It is not obvious how one could
formulate a useful version of this which replaced the upper bound by
an asymptotic formula.

{\bf Acknowledgements.}  The author is extremely grateful to Professor
Yuri Tschinkel, who introduced him to this problem, and gave him a
full description of the universal torsor.

The work described here was carried out while
the author was visiting the Max-Planck Institute for Mathematics in
Bonn, and the American Institute of Mathematics.  The hospitality and
financial support of both institutes is gratefully acknowledged.

\section{The Universal Torsor}

Our goal in this section is to use factorization information to
analyze the equation $C(\x)=0$, introducing further variables which
will be of smaller size than the original variables $x_1,\ldots,x_4$,
and which will satisfy additional equations.  Although we shall not
make any use of the fact, we note that 
these new variables describe the `Universal
Torsor' for the Cayley cubic.  For the purposes of this analysis it
will be convenient to introduce the convention that the letters
$i,j,k,l$ will denote generic distinct indices from the set $\{1,2,3,4\}$.

It is useful to begin by observing that none of the variables $x_i$
can vanish.  For if $x_i=0$, then the equation $C(\x)=0$
implies that $x_jx_kx_l=0$, so that the point $\x$ must lie on one of 
the excluded lines $X_i=X_j=0$.  We now set
\begin{equation}\label{2.1}
y_i=\hcf(x_j,x_k,x_l),
\end{equation}
The requirement that $\x$ is primitive is then equivalent to the condition
\begin{equation}\label{2.2}
\hcf(y_i,y_j)=1,
\end{equation}
According to our convention this should be taken to mean that $y_i$
and $y_j$ are coprime whenever $i$ and $j$ are distinct.
Since $y_j,y_k,y_l$ are pairwise coprime and all divide $x_i$, their
product divides $x_i$, and similarly for the other indices.  We may
therefore set
\[x_i=y_jy_ky_lz_i.\]
The definition (\ref{2.1}) now reduces to 
\[
\hcf(y_ky_lz_j,y_jy_lz_k,y_jy_kz_l)=1,
\]
In view of (\ref{2.2}), this is equivalent to the two conditions
\begin{equation}\label{2.3}
\hcf(y_i,z_i)=1,
\end{equation}
and
\begin{equation}\label{2.4}
\hcf(z_i,z_j,z_k)=1.
\end{equation}
Moreover the equation $C(\x)=0$ becomes
\begin{equation}\label{2.5}
z_2z_3z_4y_1+z_1z_3z_4y_2+z_1z_2z_4y_3+z_1z_2z_3y_4=0
\end{equation}
on recalling that none of $x_1,\ldots,x_4$ can vanish.

Our problem is therefore reduced to counting solutions of the equation
(\ref{2.5}), lying in the region
\[|y_jy_ky_lz_i|\le B,\]
and subject to the constraints (\ref{2.2}), (\ref{2.3}) and (\ref{2.4}).
Moreover solutions in which any of the variables is zero are to be
discounted, since they produce points $\x$ on one of the lines in the
surface $C$.  Similarly solutions with
\[z_jz_kz_ly_i+z_iz_kz_ly_j=0\]
are to be discounted.

We now perform a second reduction.  We begin by defining
\begin{equation}\label{2.6}
z_{ij}=z_{ji}=\hcf(z_i,z_j).
\end{equation}
In view of (\ref{2.4}) we have
\[\hcf(z_{ij},z_{ik})=1.\]
Since $z_{ij},z_{ik},z_{il}$ all divide $z_i$, and are coprime in
pairs, it follows that their product divides $z_i$.  We may therefore write
\[z_i=B_iw_i\]
where
\begin{equation}\label{2.7}
B_i=z_{ij}z_{ik}z_{il}.
\end{equation}
The definition (\ref{2.6}) then reduces to
\[\hcf(z_{ik}z_{il}w_i,z_{jk}z_{jl}w_j)=1,\]
or equivalently
\begin{equation}\label{2.8}
\hcf(w_i,w_j)=1,
\end{equation}
\begin{equation}\label{2.9}
\hcf(w_i,z_{jk})=1,
\end{equation}
and
\begin{equation}\label{2.10}
\hcf(z_{ab},z_{cd})=1,\;\;\;(\{a,b\},\,\{c,d\}\;\mbox{distinct}).
\end{equation}
The equation (\ref{2.5}) now becomes
\begin{equation}\label{2.11}
A_1w_2w_3w_3y_1+A_2w_1w_3w_4y_2+A_3w_1w_2w_4y_3+A_4w_1w_2w_3y_4=0,
\end{equation}
where
\begin{equation}\label{2.12}
A_i=z_{jk}z_{jl}z_{kl}.
\end{equation}
We therefore see that
$w_i|A_iw_jw_kw_ly_i$.  In view of (\ref{2.8}) and (\ref{2.9}) this
imples that $w_i|y_i$.  Since $w_i|z_i$ we conclude from (\ref{2.3}) that
$w_i=\pm 1$.

We now have
\[x_i=B_iy_jy_ky_lw_i,\]
with $w_i=\pm 1$. However,
in making the definitions (\ref{2.1}) and (\ref{2.6}), the highest 
common factors are
only defined up to sign.  Let us assume, temporarily, that we chose 
the variables $y_i$ and $z_{ij}$ to be positive.  We proceed to replace each
$y_i$ by $w_iy_i$, whence
\[x_i=\ep B_iy_jy_ky_l,\]
with $\ep=w_1w_2w_3w_4=\pm 1$.  Thus, if we replace $\x$ by $\ep\x$, we
obtain
\begin{equation}\label{2.13}
x_i=B_iy_jy_ky_l,
\end{equation}
where the variables $z_{ij}$ are positive but $y_i$ may be of either
sign.
After these changes the equation (\ref{2.11}) reduces to
\begin{equation}\label{2.14}
A_1y_1+A_2y_2+A_3y_3+A_4y_4=0.
\end{equation}
Moreover the condition (\ref{2.4}) is implied by (\ref{2.10}), 
while (\ref{2.3})
is equivalent to
\begin{equation}\label{2.15}
\hcf(y_i,z_{ij})=1.
\end{equation}
We may therefore summarize our conclusions as follows.
\begin{lemma}
Let $\x\in U$ be a primitive integral solution of $C(\x)=0$.  Then
either $\x$ or $-\x$ takes the form
(\ref{2.13}), with non-zero integer variables $y_i$ and positive integer
variables $z_{ij}$ constrained by the conditions (\ref{2.2}), (\ref{2.10})
and 
(\ref{2.15}), and
satisfying the equation (\ref{2.14}).  Moreover none of
$A_1y_1+A_2y_2$,
$A_1y_1+A_3y_3$ or
$A_1y_1+A_4y_4$ may vanish.

Conversely, if $y_i$ and $z_{ij}$ are as above, then the vector $\x$
given by (\ref{2.13}) will be a primitive integral solution of $C(\x)=0$
lying in $U$.
\end{lemma}

To proceed further, we note that the equation (\ref{2.14}) implies that
\[z_{ij}|z_{kl}(z_{ik}z_{il}y_j+z_{jk}z_{jl}y_i),\]
whence (\ref{2.10}) yields
\[z_{ij}|z_{ik}z_{il}y_j+z_{jk}z_{jl}y_i.\]
We therefore write
\begin{equation}\label{2.16}
z_{ik}z_{il}y_j+z_{jk}z_{jl}y_i=z_{ij}v_{ij},
\end{equation}
so that equation (\ref{2.14}) is equivalent to each of the relations
\[v_{ij}+v_{kl}=0.\]
Note that $v_{ij}=v_{ji}$, since $z_{ij}$ is also
symmetric in the indices $ij$. 
We now calculate that
\begin{eqnarray*}
z_{ij}v_{ij}z_{ik}v_{ik}&=&(z_{ik}z_{il}y_j+z_{jk}z_{jl}y_i)
(z_{ij}z_{il}y_k+z_{jk}z_{kl}y_i)\\
&=&z_{il}^2z_{ij}z_{ik}y_jy_k+z_{jk}y_i(z_{ik}z_{il}z_{kl}y_j+
z_{jl}z_{ij}z_{il}y_k+z_{jk}z_{jl}z_{kl}y_i),
\end{eqnarray*}
whence (\ref{2.14}) yields
\[z_{ij}v_{ij}z_{ik}v_{ik}=z_{il}^2z_{ij}z_{ik}y_jy_k-
z_{jk}y_i.z_{ij}z_{ik}z_{jk}y_l.\]
We therefore conclude that
\begin{equation}\label{2.18}
v_{ij}v_{ik}=z_{il}^2y_jy_k-z_{jk}^2y_iy_l.
\end{equation}

\section{The Lower Bound}

To tackle the lower bound problem in our theorem
we begin by considering solutions for which the variables $z_{ij}$ are
fixed, and relatively small, while the variables $y_i$ are
comparatively large, and lie in the dyadic ranges 
\begin{equation}\label{3.1}
Y_i< |y_i|\le 2Y_i.
\end{equation}
In the notation given by (\ref{2.7}) and (\ref{2.12}) we 
observe that the condition $\max|x_i|\le B$ is equivalent to
\begin{equation}\label{3.2}
A_iA_jA_k|y_iy_jy_k|\le BP,
\end{equation}
where
\begin{equation}\label{3.3}
P=z_{12}z_{13}z_{14}z_{23}z_{24}z_{34}=A_iB_i.  
\end{equation}
We shall choose
\begin{equation}\label{3.4}
Y_i=[\frac{(BP)^{1/3}}{2A_i}],
\end{equation}
We will then have (\ref{3.2}) whenever the $y_i$ lie in the 
ranges (\ref{3.1}).  We
shall assume moreover that
\begin{equation}\label{3.5}
P\le B^{\delta}.
\end{equation}
where $\delta$ is a small positive constant to be specified later, see
(\ref{3.14}).  We shall write 
\[\cl{N}=\cl{N}(z_{12},z_{13},z_{14},z_{23},z_{24},z_{34})  \]
for the number of solutions $(y_1,y_2,y_3,y_4)$ of (\ref{2.14}), in
the ranges (\ref{3.1}), subject to the constraints (\ref{2.2}) 
and (\ref{2.15}),
and not on any of
the lines $A_1y_i+A_iy_i=0$.

The main difficulty in establishing our lower bound comes from the
coprimality conditions (\ref{2.2}) and (\ref{2.15}).  To handle these we
begin 
by setting
\[Q=P\prod_{p\le\sqrt{\log B}}p\]
and writing $\cl{N}_1$ for the number of solutions in which (\ref{2.2}) is
replaced by the weaker condition 
\[(y_i,y_j,Q)=1,\;\;\;\mbox{for all}\;\;i\not=j.\]
We take $\cl{N}_2$ to be the number of solutions in which some
pair $y_i,y_j$ has a prime factor $p|y_i,y_j$ with $p\nmid Q$.
Clearly we then have
\begin{equation}\label{3.6}
\cl{N}\ge\cl{N}_1-\cl{N}_2.
\end{equation}

We begin by estimating $\cl{N}_1$, and first note that there can be at
most 
\begin{equation}\label{3.7}
Y_1Y_j\ll (BP)^{2/3}\ll B^{2/3+2\delta/3}
\end{equation}
solutions on one of the lines $A_1y_1+A_iy_i=A_jy_j+A_ky_k=0$, by (\ref{3.4})
and (\ref{3.5}).  This bound will turn out to be of negligible size.
We can therefore ignore the condition that solutions may not
lie on such a line.  We now proceed by picking out the coprimality
conditions with the M\"{o}bius function.  Let
\[\cl{N}_3=\cl{N}_3(d_1,\ldots,d_4;d_{12},\ldots,d_{34})\]
denote the number of solutions of the equation (\ref{2.14}), with $y_i$ in the
ranges (\ref{3.1}), and such that
$d_i|y_i$ and $d_{ij}|y_i,y_j$ for every choice of indices.  Then
\begin{equation}\label{3.8}
\cl{N}_1=\sum_{d_i|B_i}\mu(d_1)\ldots\mu(d_4)\sum_{d_{ij}|Q}
\mu(d_{12})\ldots\mu(d_{34})\cl{N}_3,
\end{equation}
since the condition (\ref{2.15}) is equivalent to $\hcf (y_i,B_i)=1$.
We must now estimate $\cl{N}_3$.
We shall write $h_1$ for the lowest common multiple of
$d_{12},d_{13},d_{14}$ and $d_1$,  and similarly for $h_2,h_3$ and $h_4$.
We may then re-interpet $\cl{N}_3$ as the number of integer triples 
$(n_1,n_2,n_3)$ for which $A_ih_i|n_i$ and $A_4h_4|n_1+n_2+n_3$, and
which lie in the region
\begin{equation}\label{3.9}
\cl{R}:\; A_iY_i<|n_i|\le 2A_iY_i\;\;\;(i\le 3),\;\;\;
A_4Y_4<|n_1+n_2+n_3|\le 2A_4Y_4.
\end{equation}
The divisibility conditions define an integer sublattice
$\ssl\le\Z^3$, such that
\[\cl{N}_3=\#(\ssl\cap\cl{R}).\]
We shall need to compute the determinant of $\ssl$, or, what is the same
thing, its index in $\Z^3$.  This is most easily done locally.  Write
\[A_ih_i=\prod_p p^{\nu(p,i)}\]
and let $\ssl_p\le\Z^3$ be the lattice for which $p^{\nu(p,i)}|n_i$
and $p^{\nu(p,4)}|n_1+n_2+n_3$.  Let $\nu(p,0)=\min_i\nu(p,i)$ and let
$\ssl_p^{(0)}\le\Z^3$ 
be defined by the conditions $p^{\nu(p,i)-\nu(p,0)}|n_i$
and $p^{\nu(p,4)-\nu(p,0)}|n_1+n_2+n_3$.  Then a moment's thought
reveals that
\[\det(\ssl_p^{(0)})=p^{\{\nu(p,1)-\nu(p,0)\}+\{\nu(p,2)-\nu(p,0)\}
+\{\nu(p,3)-\nu(p,0)\}+\{\nu(p,4)-\nu(p,0)\}},\]
and
\[\det(\ssl_p)=p^{3\nu(p,0)}\det(\ssl_p^{(0)}),\]
whence
\[\det(\ssl_p)=\frac{\prod_i p^{\nu(p,i)}}{\hcf (p^{\nu(p,i)})}.\]
We now observe that $\ssl$ is the 
intersection of the various $\ssl_p$, which
have pairwise coprime indices in $\Z^3$.  It therefore follows that
\[\det(\ssl)=\prod_p\frac{\prod_i p^{\nu(p,i)}}{\hcf (p^{\nu(p,i)})}=
\frac{\prod_i A_ih_i}{\hcf (A_ih_i)}.\]

We may choose a basis $\b{b}_1,\b{b}_2,\b{b}_3$ of $\ssl$ with
$|\b{b}_i|\ll\det(\ssl)$.  Taking $M$ to be the $3\times 3$ integer
matrix formed from the vectors $\b{b}_i$ we see that
$\ssl=M\Z^3$, and that $\det(M)=\det(\ssl)$.  If $\cl{R}$ is the
region (\ref{3.9}) then 
\[\#(\ssl\cap\cl{R})=\#(\Z^3\cap M^{-1}\cl{R}).\]
However $M^{-1}\cl{R}$ has volume ${\rm meas}(\cl{R})/\det(M)$, is
bounded by $O(1)$ planar sides, and lies in a sphere of radius $r$,
say, where
$r\ll ||M^{-1}||\max Y_i$.  Here $||M^{-1}||$ is the modulus of the
largest entry in $M^{-1}$, so that $||M^{-1}||\ll \det(\ssl)$.
It follows that
\[\cl{N}_3=\#(\ssl\cap\cl{R})=\#(\Z^3\cap M^{-1}\cl{R})=
\frac{{\rm meas}(\cl{R})}{\det(M)}+O(r^2).\]
Since $d_i|B_i$ and $d_{ij}|Q$ we have $A_ih_i\le
A_iB_iQ=PQ$, by (\ref{3.3}).  Thus (\ref{3.5}) yields
$\det(\ssl)\le P^4Q^4\ll
B^{12\delta}$, since 
\[Q\ll P\exp(O(\sqrt{\log B}))\ll B^{2\delta}.\]
We therefore deduce that
\[\cl{N}_3={\rm meas}(\cl{R})\frac{\hcf (A_ih_i)}{\prod_i A_ih_i}
+O(B^{2/3+25\delta}).\]
We now insert this into (\ref{3.8}), so that
\begin{eqnarray}\label{3.10}
\cl{N}_1&=&{\rm meas}(\cl{R})\sum_{d_i,d_{ij}}
\mu(d_1)\ldots\mu(d_4)
\mu(d_{12})\ldots\mu(d_{34})\frac{\hcf (A_ih_i)}{\prod_i A_ih_i}\nonumber\\
&&\hspace{3cm}+O(B^{2/3+26\delta}),
\end{eqnarray}
since the usual estimate for the divisor function shows that there are
$O(B^{\delta})$ divisors $d_i,d_{ij}$ in total.

It remains to consider the sum
\[\sum_{d_i,d_{ij}}\mu(d_1)\ldots\mu(d_4)
\mu(d_{12})\ldots\mu(d_{34})\frac{\hcf (A_ih_i)}{\prod_i A_ih_i}.\]
By multiplicativity we see that this is a product of local factors $e_p$, say.
For primes $p\nmid P$ we define the integer $N$, temporarily, as the number of
quadruples $(x_1,x_2,x_3,x_4)\;\mod{p}$ satisfying $\hcf(x_i,x_j,p)=1$
for $i\not=j$, and such that $x_1+x_2+x_3+x_4\equiv 0\;\mod{p}$.  We
then find, again using the M\"{o}bius function, that $N=p^3e_p$.  An 
easy computation then yields
\begin{equation}\label{3.11}
e_p=1-\frac{6}{p^2}+\frac{5}{p^3}.
\end{equation}
For the remaining primes $p$ we note that $p$ will divide exactly one
$z_{ij}$, by (\ref{2.10}), 
and we suppose without loss of generality that $z_{12}$ 
contains $p$ with exponent $e\ge
1$, say.   We then let $A'_i=A_i$ if $i=1$ or $2$, and $A'_i=p^{1-e}A_i$
for $i=3$ or $4$.  Since $p^2\nmid A_ih_i$ for $i=1$ or $2$ we have
\[\hcf (A_ih_i)=\hcf (A'_ih_i),\]
and hence
\[\frac{\hcf (A_ih_i)}{\prod_i A_ih_i}=
p^{2-2e}\frac{\hcf (A'_ih_i)}{\prod_i A'_ih_i}.\]
Now let $N$ denote, temporarily, the number of quadruples
$(x_1,x_2,x_3,x_4)\in \N^4$ satisfying the conditions
\[x_1,x_2\le p^2,\;\;\; x_3,x_4\le p,\]
\[\hcf(x_i,x_j,p)=1\;\mbox{for}\;
i\not=j,\;\;\;\hcf(x_1,p)=\hcf(x_2,p)=1,\]
and
\[x_1+x_2+px_3+px_4\equiv 0\mod{p^2}.\]
We then find, using the M\"{o}bius function once more, that
$N=p^{2e+4}e_p$, and
another easy computation then produces
\begin{equation}\label{3.12}
e_p=(1-\frac{1}{p})(1-\frac{1}{p^2})p^{-2e}.
\end{equation}

The formulae (\ref{3.11}) and (\ref{3.12}) show that
\[\sum_{d_i,d_{ij}}\mu(d_1)\ldots\mu(d_4)
\mu(d_{12})\ldots\mu(d_{34})\frac{\hcf (A_ih_i)}{\prod_i A_ih_i}\gg
P^{-2}\frac{\phi(P)}{P},\]
and since we clearly have ${\rm meas}(\cl{R})\gg BP$, from 
(\ref{3.4}) and (\ref{3.9}), we
deduce from (\ref{3.5}), (\ref{3.7}) and (\ref{3.10}) that
\begin{equation}\label{3.13}
\cl{N}_1\gg \frac{B}{P}\frac{\phi(P)}{P}
\end{equation}
providing that we take
\begin{equation}\label{3.14}
\delta=\frac{1}{84}.
\end{equation}

We turn now to $\cl{N}_2$, which we must estimate from above.  We
start by considering the
contribution from solutions in which $p|y_1,y_2$, say, with $p\nmid Q$
and $R<p\le
2R$.  We begin with the following preliminary observations. Clearly 
there are no solutions
with $R\gg Y_1$, and so we may suppose that $R\ll Y_1$.  Moreover, if
$\delta\le 1/7$ then (\ref{3.3}), (\ref{3.4})
and (\ref{3.5}) yield
\[z_{12}\le\frac{P}{A_1}\le\frac{(BP)^{1/6}}{A_1}\ll Y_1^{1/2}.\]
It therefore follows that
\[1\ll (\frac{Y_1}{R_1})^{1/2}\frac{Y_1^{1/2}}{z_{12}}=
\frac{Y_1}{R^{1/2}z_{12}}.\]
Thus
\begin{equation}\label{3.15}
1+\frac{Y_1}{Rz_{12}}\ll\frac{Y_1}{R^{1/2}z_{12}}.
\end{equation}

Since we are seeking an upper bound for $\cl{N}_2$, the coprimality
conditions can be dropped.  If we set $y_i=pt_1,y_2=pt_2$ then we have
$Y_i/R\ll|t_i|\ll Y_i/R$ for $i=1,2$.  Moreover, since $z_{12}|A_3,A_4$
we have $A_1t_1\equiv -A_2t_2\;\mod{z_{12}}$.  However $z_{12}$ is
coprime to $A_1$, so that each admissable value of $t_2$ determines
$O(1+Y_1/(Rz_{12}))$ values of $t_1$.  It therefore follows from (\ref{3.15})
that there are $O(Y_1Y_2R^{-3/2}z_{12}^{-1})$ possible pairs
$t_1,t_2$.

For each such pair we now estimate how many triples $p,y_3,y_4$ there
might be.  We put $A_1t_1+A_2t_2=z_{12}s$.  Then 
\[ps+\frac{A_3}{z_{12}}y_3+\frac{A_4}{z_{12}}y_4=0,\]
whence
\[ps+\frac{A_3}{z_{12}}y_3\equiv 0\mod{\frac{A_4}{z_{12}}}.\]
Since $A_3/z_{12}$ and $A_4/z_{12}$ are coprime, by (\ref{2.10}), 
we see that each value
of $p$ determines $y_3$ modulo $A_4/z_{12}$, producing $O(1+Y_3z_{12}/A_4)$
values.  However if $\delta\le 1/7$ then (\ref{3.4}) and (\ref{3.5}) 
suffice to show that $1\ll Y_3/A_4$, so we
will have $O(RY_3z_{12}/A_4)$ possible pairs $p,y_3$, each of which
determines at most one admissable $y_4$.

These bounds show that the range $R<p\le 2R$ contributes
\[\frac{Y_1Y_2}{R^{3/2}z_{12}}\frac{RY_3z_{12}}{A_4}=
\frac{A_1Y_1.A_2Y_2.A_3Y_3}{P^2R^{1/2}}\ll\frac{B}{PR^{1/2}}\]
to $\cl{N}_2$.  If we now sum $R\gg\sqrt{\log B}$ over powers of two
we deduce that
\[\cl{N}_2\ll \frac{B}{P}(\log B)^{-1/4}.\]
Since $\phi(P)/P\gg (\log\log P)^{-1}$ we deduce from (\ref{3.13}) that
$\cl{N}_2=o(\cl{N}_1)$ and hence, via (\ref{3.6}), that
\[\cl{N}\gg \frac{B}{P}\frac{\phi(P)}{P}.\]

We summarize our conclusions thus far as follows.
\begin{lemma}
For a given admissible set of values $z_{12},\ldots,z_{34}$ satisfying
\[P=\prod z_{ij}\le B^{1/84}\]
there are
\[\gg \frac{B}{P}\frac{\phi(P)}{P}\]
corresponding values of $y_1,y_2,y_3,y_4$.
\end{lemma}

To complete the proof of the lower bound part of our theorem, we
observe that any square-free value of $P$ will factorize into values
$z_{12},\ldots,z_{34}$ satisfying (\ref{2.10}) and
(\ref{3.3}) in exactly $d_6(P)$ ways. (Here
$d_6(\ldots)$ is the generalized divisor function.)  Thus
\[N(B)\gg\sum_{P\le B^{1/84}}\mu(P)^2d_6(P)\frac{B}{P}\frac{\phi(P)}{P}\]
and a standard estimation using Perron's formula then produces the
required bound
\[N(B)\gg B(\log B)^6.\]

\section{The Upper Bound---Basic Estimates}

In contrast to the work of the previous section, in giving an upper
bound for $N(B)$ we can ignore questions of coprimality whenever we
wish to do so.  Instead our principal technical problem will be to
control precisely the number of logarithms appearing in our
estimates.  

We shall need to understand the equations (\ref{2.16}) and (\ref{2.18}), 
and our results
are summarized as follows.
\begin{lemma}
Let real numbers $K_1,\ldots,K_7>0$ be given, and let $N_1$ denote the
number of solutions $n_i\in\N$ to the equation
\[n_1n_2n_3+n_4n_5n_6=n_7n_8\;\;\;(K_i<n_i\le
2K_i,\;\;1\le i\le 7)\]
subject to the condition
\begin{equation}\label{4.1}
\hcf(n_1n_2n_3\,,\,n_4n_5n_6)=1.
\end{equation}
Then 
\begin{equation}\label{4.2}
N_1\ll K_1K_2K_3K_4K_5K_6.
\end{equation}
Similarly, if $N_2$ is the number of solutions of
\begin{equation}\label{4.3}
n_1n_2n_3=n_4n_5n_6+n_7n_8
\end{equation}
under the same conditions, then 
\[N_2\ll K_1K_2K_3K_4K_5K_6.\]
\end{lemma}

\begin{lemma}
Let real numbers $K_1,\ldots,K_7>0$ be given, and let $N_3$ denote the
number of solutions $n_i\in \N$ of the equation
\begin{equation}\label{4.5}
n_1^2n_2n_3+n_4^2n_5n_6=n_7n_8\;\;\;(K_i<n_i\le
2K_i,\;\;1\le i\le 7)
\end{equation}
subject to the condition (\ref{4.1}).
Then 
\begin{equation}\label{4.6}
N_3\ll
K_1K_2K_3K_4K_5K_6\max\{(\frac{K_1^2K_2K_3}{K_4^2K_5K_6})^{1/4}\,,\,
(\frac{K_4^2K_5K_6}{K_1^2K_2K_3})^{1/4}\}.
\end{equation}
If $N_4$ is the corresponding number of solutions for the
equation
\begin{equation}\label{4.7}
n_1^2n_2n_3-n_4^2n_5n_6=n_7n_8\;\;\;(K_i<n_i\le
2K_i,\;\;1\le i\le 7)
\end{equation}
we have
\begin{eqnarray}\label{4.8}
N_4&\ll&\{1+\frac{\log K_1K_4}{(K_2K_3K_5K_6)^{1/3}}\}
K_1K_2K_3K_4K_5K_6\nonumber\\
&&\hspace{2cm}\mbox{}\times\max\{(\frac{K_1^2K_2K_3}{K_4^2K_5K_6})^{1/4}\,,\,
(\frac{K_4^2K_5K_6}{K_1^2K_2K_3})^{1/4}\}.
\end{eqnarray}
\end{lemma}

We may think of the bound for $N_1$, for example, as describing the
number of divisors of $n_1n_2n_3+n_4n_5n_6$ which lie in specified
dyadic ranges.  Note that we do not impose a condition on the size of
$n_8$.  We may remark that in both lemmas we can use the
standard bound for the divisor function to show that each $6$-tuple
$(n_1,\ldots,n_6)$ determines $O((\max K_i)^{\ep})$ pairs of divisors
$n_7,n_8$, for any fixed $\ep>0$.  This immediately yields the bounds
\[N_1,N_2,N_3, N_4\ll (K_1K_2K_3K_4K_5K_6)^{1+\ep},\]
so that the important aspect of Lemma 3 is the removal of the exponent
$\ep$.  It would be relatively easy to replace the $\ep$ power by a
power of a logarithm, but this would be insufficient for our
purposes. In relation to Lemma 4 we conjecture that the factor
\[\max\{(\frac{K_1^2K_2K_3}{K_4^2K_5K_6})^{1/4}\,,\,
(\frac{K_4^2K_5K_6}{K_1^2K_2K_3})^{1/4}\}\]
may be removed in both cases.  However it is not possible to delete
the term
\[1+\frac{\log K_1K_4}{(K_2K_3K_5K_6)^{1/3}}\]
in our estimate for $N_4$.  Indeed, when $K_1=K_4=K_7$ and
$K_2=K_3=K_5=K_6=1/2$ we easily find that $N_4\gg K_1K_4\log(K_1K_4)$.
Thus our bounds are not as sharp as we would like, but they are
optimal in the critical case in which
$K_1^2K_2K_3$ and $K_4^2K_5K_6$ have the same order of magnitude.

Before beginning the proofs of these results we observe that the
condition (\ref{4.1}) implies that the three terms $n_1n_2n_3$,
$n_4n_5n_6$ and $n_7n_8$ are coprime in pairs.  We shall use this fact
repeatedly without further comment, in relation to both lemmas.

In this section we shall prove Lemma 3.  The treatment of Lemma 4,
which we defer to the next section, uses some of the same principles,
but is much more involved.
We begin by considering $N_1$.  By the symmetry we may
assume that
\begin{equation}\label{4.9}
K_1K_2K_3\gg K_4K_5K_6.
\end{equation}
It is then clear that $N_1=0$
unless 
\begin{equation}\label{4.10}
\frac{K_1K_2K_3}{K_7}\ll n_8\ll \frac{K_1K_2K_3}{K_7},
\end{equation}
as we shall now assume.  We write this condition as $K_8\ll n_8\ll
K_8$.  We may then suppose, by symmetry,
that $K_7\ge K_8$, whence (\ref{4.10}) implies
that $K_8\ll (K_1K_2K_3)^{1/2}$.  We then apply the following estimate.
\begin{lemma}
Let $K_1,K_2,K_3>0$ and let $q\ll (K_1K_2K_3)^{1/2}.$  Then 
for any integer
$a$ coprime to $q$, we have
\[\#\{(n_1,n_2,n_3)\in\N^3:\, K_i<n_i\le 2K_i,\, n_1n_2n_3\equiv
a\mod{q}\}\]
\[\hspace{5cm}\ll K_1K_2K_3/\phi(q).\]
\end{lemma}

We shall prove this in a moment.  However if we apply it to the current
situation we see, on taking $q=n_8$ and summing over
$n_4,n_5,n_6$ and $n_8$, that
\[N_1\ll K_1K_2K_3K_4K_5K_6\sum_{K_8\ll n_8\ll K_8}\phi(n_8)^{-1}.\]
The required bound (\ref{4.2}) now follows, since
\[\sum_{K_8\ll n_8\ll K_8}\phi(n_8)^{-1}\ll 1.\]

To handle $N_2$ we note that we automatically have (\ref{4.9}) if there are to
be any solutions.  We can then proceed exactly as before providing
that
\[K_1K_2K_3\ge 16K_4K_5K_6,\]
since this is enough to ensure 
that (\ref{4.10}) holds.  It therefore remains to
consider the case in which
\begin{equation}\label{4.12}
K_1K_2K_3 \ll K_4K_5K_6\ll K_1K_2K_3.
\end{equation}
In this case we shall assume that
\[K_1\ge K_2\ge K_3,\;\;\;\mbox{and}\;\;\;K_4\ge K_5\ge K_6,\]
as we may, by the symmetry.  It follows in particular that
\begin{equation}\label{4.13}
K_2K_3K_5K_6\le (K_1K_2K_3K_4K_5K_6)^{2/3}.
\end{equation}
We now write $N_{2,a}(q)$ for the the number of
solutions $(n_1,\ldots,n_6)$ corresponding to each value $n_7=q$, so
that
\begin{equation}\label{4.14}
N_2\ll\sum_{K_7<q\le 2K_7}N_{2,a}(q).
\end{equation}
Moreover, if we set
\begin{equation}\label{4.15}
K_8=K_1K_2K_3/K_7,
\end{equation}
then it is apparent that we must have $n_8\ll K_8$ in any solution of
(\ref{4.3}).  Thus
\begin{equation}\label{4.16}
N_2\ll\sum_{q\ll K_8}N_{2,b}(q),
\end{equation}
where $N_{2,b}(q)$ counts the solutions corresponding to a given value
$n_8=q$.
We plan to use (\ref{4.14}) when $K_7\ll (K_1K_2K_3)^{1/2}$.  
If this condition 
fails to hold we must have
$K_8\ll (K_1K_2K_3)^{1/2}$, in which case we shall employ (\ref{4.16}).

We now introduce the following result, which 
is part of Lemma 3 of the author's work \cite{sfn}.
\begin{lemma}
Let $\b{v}\in\Z^3$ be a primitive vector, and let $H_i>0$ for
$i=1,2,3$ be given.  Then the number of primitive vectors
$\b{x}\in\Z^3$ for which
$\b{v}.\b{x}=0$, and which lie in the box $|x_i|\le H_i\;(i=1,2,3)$, is at most
\[4+12\pi\frac{H_1H_2H_3}{\max H_i|v_i|}\le 4+12\pi\frac{H_1H_2}{|v_3|}.\]
\end{lemma}
Recall that an integer vector is said to be primitive if its
coordinates have no common factor.  In our applications this condition
will 
be a consequence of (\ref{2.2}), (\ref{2.10}) and (\ref{2.15}).

To bound $N_{2,a}(q)$ we write the condition (\ref{4.3}) 
as $\b{v}.\b{x}=0$ where
$\b{v}=(n_2n_3,-n_5n_6,-q)$ and $\b{x}=(n_1,n_4,n_8)$.  We set
$H_1=2K_1, H_2=2K_4$ and 
\[H_3=8K_1K_2K_3/q.\]
Then Lemma 6 produces the bound $O(1+K_1K_4/q)$ for the number of
triples $(n_1,n_4,n_8)$ and it follows on 
summing over $n_2,n_3,n_5,n_6$ and $q$ that
\begin{equation}\label{4.17}
N_2\ll\sum_{K_7<q\le 2K_7}N_{2,a}(q)\ll K_2K_3K_5K_6K_7+K_1K_2K_3K_4K_5K_6.
\end{equation}

Alternatively, we may use $N_{2,b}(q)$, and write (\ref{4.3}) as
$\b{v}.\b{x}=0$ with
$\b{v}=(n_2n_3,-n_5n_6,-q)$ and $\b{x}=(n_1,n_4,n_7)$.  We set
$H_1=2K_1, H_2=2K_4$ as before, and $H_3=2K_7$.
This time Lemma 6 produces a bound 
\[\ll 1+\frac{K_1K_4K_7}{\max H_i|v_i|}\ll 1+\frac{K_1K_4K_7}{K_1K_2K_3}\]
for the number of
triples $(n_1,n_4,n_7)$. On 
summing over $n_2,n_3,n_5,n_6$ and $q$ we then find that
\[N_2\ll\sum_{q\ll K_8}N_{2,b}(q)\ll K_2K_3K_5K_6K_8+K_4K_5K_6K_7K_8.\]
In view of (\ref{4.15}) we may combine this with (\ref{4.17}) to deduce that
\[N_2\ll K_2K_3K_5K_6\min(K_7,K_8)+K_1K_2K_3K_4K_5K_6.\]
Since we have
\[\min(K_7,K_8)\ll (K_1K_2K_3)^{1/2}\ll(K_1K_2K_3K_4K_5K_6)^{1/4},\]
by (\ref{4.12}), we then deduce from (\ref{4.13}) that
\begin{eqnarray*}
N_2&\ll &(K_1K_2K_3K_4K_5K_6)^{2/3}.(K_1K_2K_3K_4K_5K_6)^{1/4}\\
&&\hspace{3cm}\mbox{}+
K_1K_2K_3K_4K_5K_6\\
&\ll &K_1K_2K_3K_4K_5K_6,
\end{eqnarray*}
which completes the proof of our bound for $N_2$.

We must now establish Lemma 5.  To do this we refer to the author's work
\cite{d3} on the divisor function $d_3(n)$ in arithmetic
progressions.  If we write
\[N(K_1,K_2,K_3;a,q)\hspace{8cm}\]
\[=\#\{(n_1,n_2,n_3)\in\N^3:\, K_i<n_i\le 2K_i,\, n_1n_2n_3\equiv
a\mod{q}\}\]
then the analysis of \cite[\S 7]{d3} suffices to show that
\[N(K_1,K_2,K_3;a,q)=C(K_1,K_2,K_3;q)+O(E),\]
with $C(K_1,K_2,K_3;q)$ independent of $a$, and an error term
\[E=
(K_1K_2K_3)^{3/4+\ep}q^{-13/24}+(K_1K_2K_3)^{7/9+\ep}q^{-7/12}+
(K_1K_2K_3)^{10/13+\ep}q^{-15/26}\]
\[\hspace{3cm}+(K_1K_2K_3)^{46/57+\ep}q^{-12/19}
+(K_1K_2K_3)^{86/107+\ep}q^{-66/107},\]
for any fixed $\ep>0$.  Since $q\ll(K_1K_2K_3)^{1/2}$ we deduce that
\[N(K_1,K_2,K_3;a,q)
=C(K_1,K_2,K_3;q)+O(K_1K_2K_3q^{-1}).\]
We may now average over $a$ coprime to $q$ to find that
\[\phi(q)^{-1}\#\{(n_1,n_2,n_3)\in\N^3:\, K_i<n_i\le 2K_i,\,
\hcf(n_1n_2n_3,q)=1\}\]
\[\hspace{5cm}=C(K_1,K_2,K_3;q)+O(K_1K_2K_3q^{-1}).\]
We deduce that $C(K_1,K_2,K_3;q)\ll K_1K_2K_3/\phi(q)$, and Lemma
5 follows.  The reader should note that the work of
Friedlander and Iwaniec \cite{FI} could have been used equally
effectively at this point.

\section{The Proof of Lemma 4}

By symmetry, we may suppose at the outset that
\begin{equation}\label{4.18}
K_1^2K_2K_3\gg K_4^2K_5K_5.
\end{equation}
We shall write
$N_{3,a}(q)$ for the number of
solutions $(n_1,\ldots,n_6)$ corresponding to each value $n_7=q$, so
that
\begin{equation}\label{4.19}
N_3\ll\sum_{K_7<q\le 2K_7}N_{3,a}(q).
\end{equation}
Moreover, if we set
\[K_8=K_1^2K_2K_3/K_7,\]
then it is apparent that we must have $K_8\ll n_8\ll K_8$ in any solution of
(\ref{4.5}).  Thus
\begin{equation}\label{4.20}
N_3\ll\sum_{K_8\ll q\ll K_8}N_{3,b}(q),
\end{equation}
where $N_{3,a}(q)$ counts the solutions corresponding to a given value
$n_8=q$.

We plan to use (\ref{4.19}) when $K_7\ll
K_1(K_2K_3)^{1/2}$.  If this condition fails to hold we must have
$K_8\ll K_1(K_2K_3)^{1/2}$, in which case we shall employ (\ref{4.20}).

To bound $N_{3,a}(q)$ we write the condition (\ref{4.5}) 
as $\b{v}.\b{x}=0$ where
$\b{v}=(n_1^2n_2,n_4^2n_5,-q)$ and $\b{x}=(n_3,n_4,n_8)$, say.  We set
$H_1=2K_3, H_2=2K_6$ and 
\[H_3=32K_1^2K_2K_3/q.\]
Then Lemma 6 produces the bound $O(1+K_3K_6/q)$ for the number of
triples $(n_3,n_6,n_8)$ and it follows on 
summing over $n_1,n_2,n_4$ and $n_5$ that
\begin{equation}\label{4.21}
N_{3,a}(q)\ll K_1K_2K_4K_5+K_1K_2K_3K_4K_5K_6/q.
\end{equation}
In a precisely analogous way we find that
\begin{equation}\label{4.22}
N_{3,b}(q)\ll K_1K_3K_4K_6+K_1K_2K_3K_4K_5K_6/q.
\end{equation}
We may also use a vector $\b{x}$ involving $n_1$ and $n_4$.  To do
this, we let $t\in[0,q)$ run over the solutions of the quadratic congruence
\[t^2n_2n_3+n_5n_6\equiv 0\mod{q},\]
and we write $\rho(q;n_2n_3,n_5n_6)$ for the number of such solutions $t$.
We then see that, for fixed $n_2,n_3,n_5,n_6$ and $q$, we must have
$n_1\equiv tn_4\;\mod{q}$ for some value of $t$.  This leads to an equation
$\b{v}.\b{x}=0$ with $\b{v}=(1,-t,q)$ and $\b{x}=(n_1,n_4,m)$, with
size restrictions given by $H_1=2K_1,H_2=2K_4$ and $H_3=K_1+K_4$,
say.  Thus Lemma 6 produces a bound $O(1+K_1K_4/q)$ for the number of
solutions $n_1,n_4$ corresponding to a given value of $t$.  We
therefore obtain an estimate
\begin{equation}\label{4.23}
N_{3,a}(q)\ll(1+K_1K_4/q)\sum_{n_2,n_3,n_5,n_6}\rho(q;n_2n_3,n_5n_6).
\end{equation}

Our next task is evidently to examine averages of the function $\rho$.
Suppose that $\hcf (a,b)=1$.  Then for odd $q$ we have
\[\rho(q;a,b)=\sum_{d|q}\mu(d)^2(\frac{-ab}{d}),\]
where $(ab/d)$ is the Jacobi symbol.  We then see that
\[\rho(q;a,b)\le 4\sum_{d|q}\mu(d)^2(\frac{-ab}{d})\]
whether $q$ is even or odd,
where we take the Jacobi symbol to vanish for even $d$.  We also note
that the sum on the right is non-negative
when $ab$ and
$q$ are not coprime.
Our aim is to estimate
\[\sum_{q\le Q}\sum_{n_2,n_3,n_5,n_6}\rho(q;n_2n_3,n_5n_6)=S,\]
say.
It will facilitate our argument to average over all $4$-tuples
$(n_2,n_3,n_5,n_6)$ in the relevant ranges, and not just those
satisfying the coprimality condition (\ref{4.1}).
In view of the above remarks we clearly have
\begin{eqnarray}\label{4.24}
S&\ll&\sum_{e\le Q}\sum_{d\le Q/e}\mu(d)^2\sum_{n_2,n_3,n_5,n_6}
(\frac{-n_2n_3n_5n_6}{d})\nonumber\\
&\ll& KQ+\sum_{e\le Q}S(e),
\end{eqnarray}
where
\[K=K_2K_3K_5K_6\]
and
\[S(e)=\sum_{1\not=d\le Q/e}\mu(d)^2\sum_{n_2,n_3,n_5,n_6}
(\frac{-n_2n_3n_5n_6}{d}).\]
An immediate application of the author's large sieve inequality for
real character sums \cite[Corollary 4]{rc} shows that
\begin{equation}\label{4.25}
S(e)\ll (KQ/e)^{\ep}\{K(Q/e)^{1/2}+K^{1/2}(Q/e)\}
\end{equation}
for any fixed $\ep>0$.
If we use the
P\'{o}lya-Vinogradov inequality, we find that
\begin{eqnarray*}
S(e)&=&\sum_{1\not=d\le Q/e}\mu(d)^2\sum_{n_2,n_3,n_5}
(\frac{-n_2n_3n_5}{d})\sum_{n_6}(\frac{n_6}{d})\\
&\ll& K_2K_3K_5\sum_{1\not=d\le Q/e}d^{1/2+\ep}\\
&\ll& K_2K_3K_5(Q/e)^{3/2+\ep}.
\end{eqnarray*}
In the same way we find that
\[S(e)\ll\frac{K_2K_3K_5K_6}{K_i}(Q/e)^{3/2+\ep}\]
for any index $i=2,3,5,6$.  It therefore follows on taking $K_i$ as
the maximum of $K_2,K_3,K_5$ and $K_6$, that
\begin{equation}\label{4.26}
S(e)\ll K^{3/4}(Q/e)^{3/2+\ep}.
\end{equation}
Alternatively, if $N$ is not a square, 
we may use the P\'{o}lya-Vinogradov to derive the bound
\begin{eqnarray}\label{4.27}
\sum_{d\le Q/e}\mu(d)^2(\frac{N}{d})&=&
\sum_{d\le Q/e}\sum_{h^2|d}\mu(h)(\frac{N}{d})\nonumber\\
&=&\sum_{h\le(Q/e)^{1/2}}\mu(h)(\frac{N}{h^2})\sum_{k\le Q/eh^2}
(\frac{N}{k})\nonumber\\
&\ll&\sum_{h\le(Q/e)^{1/2}}N^{1/2}\log N\nonumber\\
&\ll& (QN/e)^{1/2}\log N.
\end{eqnarray}
We can use this estimate to find that
\begin{equation}\label{4.28}
S(e)\ll K^{3/2+\ep}(Q/e)^{1/2},
\end{equation}
since $-n_2n_3n_5n_6$ is never a square.
Comparing this bound with (\ref{4.25}) and (\ref{4.26}) we find that
\[S(e) \ll (KQ/e)^{\ep}\min,\]
where
\begin{eqnarray*}
\min&=&\min\{K(Q/e)^{1/2}+K^{1/2}(Q/e)\,,\,
K^{3/4}(Q/e)^{3/2}\,,\,K^{3/2}(Q/e)^{1/2}\}\\
&\ll&\min\{K(Q/e)^{1/2}\,,\,K^{3/4}(Q/e)^{3/2}\,,\,
K^{3/2}(Q/e)^{1/2}\}\\
&&\hspace{1cm}+\min\{K^{1/2}(Q/e)\,,\,
K^{3/4}(Q/e)^{3/2}\,,\,K^{3/2}(Q/e)^{1/2}\}\\
&\ll&\{K(Q/e)^{1/2}\}^{3/5}\{K^{3/4}(Q/e)^{3/2}\}^{2/5}\\
&&\hspace{1cm}+\{K^{1/2}(Q/e)\}^{2/3}\{K^{3/2}(Q/e)^{1/2}\}^{1/3}\\
&\ll&(KQ/e)^{9/10}+(KQ/e)^{5/6}\\
&\ll&(KQ/e)^{9/10}.
\end{eqnarray*}
It follows that
\begin{equation}\label{4.29}
S(e)\ll (KQ/e)^{10/11},
\end{equation}
say.  Finally we insert this into (\ref{4.24}) to deduce that
\begin{equation}\label{4.30}
S\ll KQ+K^{10/11}Q^{10/11}\sum_{e\le Q}e^{-10/11}\ll
KQ+K^{10/11}Q\ll KQ.
\end{equation}

The above bound allows us to conclude from (\ref{4.23}) that
\[\sum_{Q/2< q\le Q}N_{3,a}(q)\ll
(1+K_1K_4/Q)KQ=K_2K_3K_5K_6Q+K_1K_2K_3K_4K_5K_6.\]
On the other hand, (\ref{4.21}) and (\ref{4.22}) yield
\[\sum_{Q/2< q\le Q}N_{3,a}(q)\ll K_1K_2K_4K_5Q+K_1K_2K_3K_4K_5K_6\]
and
\[\sum_{Q/2< q\le Q}N_{3,a}(q)\ll K_1K_3K_4K_6Q+K_1K_2K_3K_4K_5K_6.\]
Taking the minimum of these, and assuming that $Q\ll K_1(K_2K_3)^{1/2}$,
we obtain an estimate
\begin{eqnarray}\label{4.31}
\sum_{Q/2< q\le Q}N_{3,a}(q)&\ll& 
Q\min\{K_2K_3K_5K_6\,,\,K_1K_2K_4K_5\,,\,K_1K_3K_4K_6\}\nonumber\\
&&\hspace{3cm}\mbox{}+
K_1K_2K_3K_4K_5K_6\nonumber\\
&\ll&Q\{K_2K_3K_5K_6\}^{1/2}\{K_1K_2K_4K_5\}^{1/4}\{K_1K_3K_4K_6\}^{1/4}
\nonumber\\
&&\hspace{3cm}\mbox{}+K_1K_2K_3K_4K_5K_6\nonumber\\
&\ll& \rule{0mm}{6mm}QK_1^{1/2}K_2^{3/4}K_3^{3/4}K_4^{1/2}K_5^{3/4}K_6^{3/4}+
K_1K_2K_3K_4K_5K_6\nonumber\\
&\ll& \rule{0mm}{6mm}
K_1(K_2K_3)^{1/2}.K_1^{1/2}K_2^{3/4}K_3^{3/4}K_4^{1/2}K_5^{3/4}K_6^{3/4}
\nonumber\\
&&\hspace{3cm}\mbox{}+
K_1K_2K_3K_4K_5K_6\nonumber\\
&\ll&\rule{0mm}{6mm}K_1^{3/2}K_2^{5/4}K_3^{5/4}K_4^{1/2}K_5^{3/4}K_6^{3/4}
\end{eqnarray}
in view of our assumption (\ref{4.18}).  This gives a satisfactory 
bound for (\ref{4.19}).
We may handle $N_{3,b}(q)$ in a precisely analogous way, 
thereby completing our treatment of (\ref{4.6}).

The equation (\ref{4.7}) introduces a couple of further difficulties.
Firstly, the bound (\ref{4.27}) is only valid when $N$ is not a square.
Previously we took $N=-n_2n_3n_5n_6$, which can never be a square.
However, if $N=n_2n_3n_5n_6$, we must allow for the case in which
$n_2n_3n_5n_6$ is a square.  The effect of this is to change the
estimate (\ref{4.28}) into
\[S(e)\ll K^{3/2+\ep}(Q/e)^{1/2}+K^{1/2+\ep}(Q/e).\]
The additional term contributes $O(K^{3/2}Q/e)$, say, 
to (\ref{4.29}), whence (\ref{4.30}) becomes
\begin{equation}\label{4.32a}
S\ll KQ\{1+\frac{\log(K_1K_2K_3K_4K_5K_6)}{K^{1/3}}\}\ll
KQ\{1+\frac{\log K_1K_4}{K^{1/3}}\}.
\end{equation}
This introduces the extra factor we see in (\ref{4.8}).

The second difficulty is that if 
\begin{equation}\label{4.32}
K_4^2K_5K_6\gg K_1^2K_2K_3\gg K_4^2K_5K_6
\end{equation}
we may no longer have the lower bound $q\gg K_8$ to use in the estimate
\begin{equation}\label{4.33a}
N_4\ll\sum_{q\ll K_8}N_{4,b}(q).
\end{equation}
We
therefore assume now that (\ref{4.32}) holds, and investigate the quantity
$N_{4,b}(q)$ further.
Since
\[n_1^2n_2n_3-n_4^2n_5n_6=n_7q\]
in this context, with $n_7\le 2K_7$, we can apply Lemma 6 with
\[\b{v}=(n_1^2n_2,-n_4^2n_5,q),\;\;\;\b{x}=(n_3,n_6, -n_7),\]
and with
$H_1=2K_3, H_2=2K_6$ and $H_3=2K_7$.  Thus there are
\[\ll 1+\frac{H_1H_2H_3}{\max H_i|v_i|}\ll
1+\frac{K_3K_6K_7}{K_1^2K_2K_3}\]
solutions $\b{x}$.  Summing over $n_1,n_2,n_4,n_5$ yields
\begin{equation}\label{4.33}
N_{4,b}(q)\ll K_1K_2K_4K_5(1+\frac{K_6K_7}{K_1^2K_2}).
\end{equation}
Similarly one can show that
\begin{equation}\label{4.34}
N_{4,b}(q)\ll K_1K_3K_4K_6(1+\frac{K_5K_7}{K_1^2K_2}).
\end{equation}

As before we need also an estimate in which we treat $n_1$ and $n_4$
as variables.  By the argument used before we can produce
$\rho(q;n_2n_3,-n_5n_6)$ congruence conditions $n_1\equiv tn_4\;\mod{q}$.
Each of these defines a lattice $\ssl\subseteq\Z^2$ of points
$(n_1,n_4)$.  Moreover we will have $\det(\ssl)=q$.  The points
$(n_1,n_4)$ satisfy $n_1\le 2K_1$ and $n_4\le 2K_4$.  Additionally we
have 
\[|1-\frac{n_4^2n_5n_6}{n_1^2n_2n_3}|\leq\frac{2K_7q}{K_2^2K_2K_3},\]
whence
\[|1-\frac{n_4\sqrt{n_5n_6}}{n_1\sqrt{n_2n_3}}|
\leq\frac{2K_7q}{K_1^2K_2K_3}.\]
In view of our assumption (\ref{4.32}) this may be written as
\[|n_1-\alpha n_4|\leq C\frac{K_7q}{K_1K_2K_3},\]
for some $\alpha=\alpha(n_2,n_3,n_4,n_5)$ and some absolute constant
$C$.  The above inequality, along with the condition $|n_4|\leq 2K_4$,
defines a parallelogram of area 
\[8C\frac{K_4K_7q}{K_1K_2K_3}=A,\]
say, centred on the origin.  This parallelogram may be mapped to a
square $S$, centred on the origin, and having the same area $A$,
by a projective mapping $M$ say, of determinant $1$. Enclose $S$ by a
disc $D$ of area $\pi A/2$, and consider the ellipse $E=M^{-1}D$.  This
also has area $\pi A/2$.  Moreover it 
contains the original parallelogram, and is
centred at the origin.  We are therefore in a position to apply the
following result, due to the author \cite[Lemma 2]{sfn}.
\begin{lemma}
Let $\ssl\subseteq\R^2$ be a lattice, and let $E$ be an ellipse,
centred on the origin, together with its interior.  Then
\[\#(\ssl\cap E)\le 4(1+\frac{{\rm meas}(E)}{\det(\ssl)}).\]
\end{lemma}
This lemma allows us to conclude that there are 
\[\ll 1+\frac{K_4K_7}{K_1K_2K_3}\]
pairs $(n_1,n_4)$ for each set of values $t,n_2,n_3,n_5,n_6,q$.
We may now procced as before, using (\ref{4.32a})
to deduce that
\begin{equation}\label{4.35}
\sum_{q\ll K_8}N_{4,b}(q)\ll \tau K_2K_3K_5K_6(1+\frac{K_4K_7}{K_1K_2K_3})K_8,
\end{equation}
where we have set
\[\tau=1+\frac{\log K_1K_4}{K^{1/3}}\]
for convenience.
We now deduce from (\ref{4.33a}),
(\ref{4.33}), (\ref{4.34}) and (\ref{4.35}) that
\begin{eqnarray*}
N_4&\ll&\sum_{q\ll K_8}N_{4,b}(q)\\
&\ll& \tau\min(K_1K_2K_4K_5K_8\,,\,K_1K_3K_4K_6K_8\,,\,K_2K_3K_5K_6K_8)\\
&&\hspace{3cm}+
\tau K_1K_2K_3K_4K_5K_6\frac{K_7K_8}{K_1^2K_2K_3}\\
&\ll& \tau K_1^{1/2}K_2^{3/4}K_3^{3/4}K_4^{1/2}K_5^{3/4}K_6^{3/4}K_8+
\tau K_1K_2K_3K_4K_5K_6
\end{eqnarray*}
as in the proof of (\ref{4.31}).  Since we only need 
(\ref{4.33a}) for the case $K_8\ll
K_1(K_2K_3)^{1/2}$, the required bound (\ref{4.8}) follows.

\section{Proof of the Upper Bound}

We shall specify dyadic ranges
\[X_i< |x_i|\le 2X_i\]
for the original variables $x_1,x_2,x_3,x_4$,
and 
\[Z_{ij}<z_{ij}\le 2Z_{ij}\]
for the variables introduced in \S 2, and we write
\[\cl{N}(X_1,\ldots,X_4;Z_{12},\dots,Z_{34})=\cl{N}\]
for the corresponding
contribution to $N(B)$.  We obviously have
\begin{equation}\label{5.1a}
X_i\ll B,
\end{equation}
and the relation (\ref{2.13}) implies that
\begin{equation}\label{5.1b}
Z_{ij}Z_{ik}Z_{il}\ll X_i.
\end{equation}

We shall find it convenient to re-order the indices so that 
\begin{equation}\label{5.1}
X_1\ge X_2\ge X_3\ge X_4.
\end{equation}
Since any solution will have
\[\frac{1}{2X_4}\le|\frac{1}{x_4}|=
|\frac{1}{x_1}+\frac{1}{x_2}+\frac{1}{x_3}|\le
|\frac{1}{x_1}|+|\frac{1}{x_2}|+|\frac{1}{x_3}|\le
\frac{1}{X_1}+\frac{1}{X_2}+\frac{1}{X_3}
\le\frac{3}{X_3}\]
we deduce that $\cl{N}=0$ unless
\begin{equation}\label{5.2}
6X_4\ge X_3\ge X_4,
\end{equation}
as we henceforth assume.  Moreover we have
\[y_i^3=\frac{x_jx_kx_lB_i}{(x_iA_i)^2}
=\frac{x_1x_2x_3x_4P}{(x_iA_i)^3},\]
with the notations (\ref{2.7}), (\ref{2.12}) and (\ref{3.3}), so that 
\[Y_i\ll |y_i|\ll Y_i\]
where 
\begin{equation}\label{5.3}
Y_i=\frac{(XF)^{1/3}}{X_iZ_{jk}Z_{jl}Z_{kl}},
\end{equation}
with
\[X=X_1X_2X_3X_4\;\;\;\mbox{and}\;\;\;
F=Z_{12}Z_{13}Z_{14}Z_{23}Z_{24}Z_{34}.\]

We begin by applying Lemma 3 to the equation (\ref{2.16}), to show that there
are
\[\ll Z_{ik}Z_{il}Z_{jk}Z_{jl}Y_iY_j\]
possible sets of values for 
$z_{ik},z_{il},z_{jk},z_{jl},z_{ij},y_i,y_j,v_{ij}$.  For each set of
values we proceed to examine (\ref{2.18}), which we write in the form
$\b{v}.\b{x}=0$ with 
\[\b{v}=(v_{ij},-z_{il}^2y_j,z_{jk}^2y_i)\;\;\;\mbox{and}
\;\;\; \b{x}=(v_{ik},y_k,y_l).\]
In view of (\ref{2.2}), (\ref{2.10}) and (\ref{2.15}) both $\b{v}$ and
$\b{x}$ will be primitive.  Moreover, the equation
\[z_{ij}z_{il}y_k+z_{jk}z_{kl}y_i=z_{ik}v_{ik},\]
which is an example of (\ref{2.16}), yields
\[z_{ik}v_{ik}\ll\max\{Z_{ij}Z_{il}Y_k\,,\,Z_{jk}Z_{kl}Y_i\}
=\frac{(XF)^{1/3}}{Z_{jl}}
\max\{\frac{1}{X_{k} }\,,\,\frac{1}{X_{i} }\}.\]
Thus Lemma 6 may be applied with 
\[H_1=c\frac{(XF)^{1/3}}{Z_{ik}Z_{jl}}
\max\{\frac{1}{X_{k} }\,,\,\frac{1}{X_{i} }\},\;\;\;
H_2=cY_k,\;\;\;
H_3=cY_l,\]
for a suitable constant $c$.
Since the remaining value $z_{kl}$ is
determined by (\ref{2.14}), there are 
\[\cl{N}\ll Z_{ik}Z_{il}Z_{jk}Z_{jl}Y_iY_j\{1+\frac{H_1H_2H_3}
{\max(H_2V_2\,,\,H_3V_3)}\}\]
solutions to (\ref{2.14}) in total, where
\[V_2=Z_{il}^2Y_j\;\;\;\mbox{and}\;\;\;
V_3=Z_{jk}^2Y_i.\]
We may then calculate, using (\ref{5.3}), that the above bound is
\[\ll \frac{(XF)^{2/3}}{Z_{kl}^2X_iX_j }+
\max(\frac{1}{X_k}\,,\,\frac{1}{X_i})\min(X_iX_l \,,\,X_jX_k ).\]
Since this estimate is valid for any choice of $i,j,k,l$ we may
interchange $i$ with $j$, and $k$ with $l$,
to deduce that 
\[\cl{N}\ll \frac{(XF)^{2/3}}{Z_{kl}^2X_iX_j }+
\max(\frac{1}{X_l}\,,\,\frac{1}{X_j})\min(X_jX_k \,,\,X_iX_l).\]
We now observe that our assumption (\ref{5.1}) 
implies that $\min(X_iX_l,X_jX_k)
\le X_2X_3$ and that either 
\[\max(\frac{1}{X_k}\,,\,\frac{1}{X_i})\le\frac{1}{X_2}\]
or
\[\max(\frac{1}{X_l}\,,\,\frac{1}{X_j})\le\frac{1}{X_2}.\]
It follows that
\[\cl{N}(X_1,\ldots,X_4;Z_{12},\dots,Z_{34})
\ll\frac{(XF)^{2/3}}{Z_{kl}^2X_iX_j }+X_3.\]
We apply this with $i=1,\, j=4,\,k=2,\,l=3$, so that
\[\cl{N}(X_1,\ldots,X_4;Z_{12},\dots,Z_{34})
\ll \frac{(XF)^{2/3}}{Z_{23}^2X_1X_4}+X_3,\]
and again with $i=2,\, j=3,\,k=1,\,l=4$, so that
\[\cl{N}(X_1,\ldots,X_4;Z_{12},\dots,Z_{34})
\ll \frac{(XF)^{2/3}}{Z_{14}^2X_2X_3}+X_3.\]
Since 
\begin{equation}\label{5.4}
\min(A,B)\le (AB)^{1/2}
\end{equation}
this yields
\begin{eqnarray}\label{5.5}
\cl{N}(X_1,\ldots,X_4;Z_{12},\dots,Z_{34})
&\ll &\min\{\frac{(XF)^{2/3}}{Z_{23}^2X_1X_4}\,,\,
\frac{(XF)^{2/3}}{Z_{14}^2X_2X_3}\}+X_3\nonumber\\
&\ll&
\frac{(XF)^{2/3}}{Z_{14}Z_{23}X^{1/2}}+X_3\nonumber\\
&=&\frac{X^{1/6}F^{2/3}}{Z_{14}Z_{23}}+X_3.
\end{eqnarray}

For an alternative estimate we begin by applying Lemma 4 to the
equation (\ref{2.18}), to show that the number of 
possible sets of values for
$v_{ij},v_{ik},z_{il},$ $z_{jk},y_i,y_j,y_k,y_l$ is
\begin{eqnarray*}
&\ll& \{1+\frac{\log X}{(Y_iY_jY_kY_l)^{1/3}}\}
Z_{il}Z_{jk}Y_iY_jY_kY_l
\max\{(\frac{Z_{il}^2Y_jY_k}{Z_{jk}^2Y_iY_l})^{1/4}\,,\,
(\frac{Z_{jk}^2Y_iY_l}{Z_{il}^2Y_jY_k})^{1/4}\}\\
&\ll & \sigma 
Z_{il}Z_{jk}X^{1/3}F^{-2/3}\max\{(\frac{X_iX_l}{X_jX_k})^{1/4}\,,\,
(\frac{X_jX_k}{X_iX_l})^{1/4}\},
\end{eqnarray*}
where
\[\sigma=1+\frac{\log X}{(XF^{-2})^{1/9}}.\]
For each such set of values we write (\ref{2.16}) in the form $\b{v}.\b{x}=0$ 
with 
\[\b{v}=(z_{il}y_j,z_{jk}y_i,-v_{ij})\;\;\;\mbox{and}\;\;\;
\b{x}=(z_{ik},z_{jl},z_{ij}).  \]
In view of (\ref{2.2}), (\ref{2.10}) and (\ref{2.15}) both $\b{v}$ and
$\b{x}$ will be primitive.  We can therefore apply Lemma 6 with
\[H_1=2Z_{ik},\;\;\;H_2=2Z_{jl},\;\;\;H_3=2Z_{ij},\]
to deduce that
there are
\begin{eqnarray*}
&\ll& 1+\frac{H_1H_2H_3}{\max(H_1|v_1|\,,\,H_2|v_2|)}\\
&\ll& 1+\frac{Z_{ik}Z_{jl}Z_{ij}}
{\max(Z_{ik}Z_{il}Y_j\,,\,Z_{jl}Z_{jk}Y_i)}\\
&\ll& 1+\frac{Z_{ik}Z_{jl}Z_{ij}Z_{kl}}{(XF)^{1/3}}\min(X_i\,,\,X_j)
\end{eqnarray*}
corresponding solutions $z_{ik},z_{jl},z_{ij}$.

We apply these estimates with $i=1,\, j=3,\,k=2,\,l=4$, so that
\[\max\{(\frac{X_iX_l}{X_jX_k})^{1/4}\,,\,
(\frac{X_jX_k}{X_iX_l})^{1/4}\}\ll (X_1/X_2)^{1/4},\]
and $\min(X_i\,,\,X_j)=X_3$, in view of (\ref{5.1}) and 
(\ref{5.2}).  Since the final remaining
value $z_{kl}=z_{24}$ is now determined by (\ref{2.14}) it follows that
\begin{eqnarray*}
\cl{N}&\ll &\sigma Z_{14}Z_{23}X^{1/3}F^{-2/3}(X_1/X_2)^{1/4}\left(
1+\frac{Z_{12}Z_{34}Z_{13}Z_{24}}{(XF)^{1/3}}X_3\right)\\
&=&\sigma Z_{14}Z_{23}X^{1/3}F^{-2/3}(X_1/X_2)^{1/4}
+\sigma (X_1/X_2)^{1/4}X_3.
\end{eqnarray*}
We now combine this with (\ref{5.5}), using the inequality (\ref{5.4}) 
again, to deduce that
\begin{eqnarray*}
\cl{N}&\ll&\sigma\min\{
\frac{X^{1/6}F^{2/3}}{Z_{14}Z_{23}}\,,\,
Z_{14}Z_{23}X^{1/3}F^{-2/3}(X_1/X_2)^{1/4}\}+\sigma
(X_1/X_2)^{1/4}X_3\\
&\ll&\sigma X^{1/4}(X_1/X_2)^{1/8}+\sigma (X_1/X_2)^{1/4}X_3.
\end{eqnarray*}

We are finally in a position to sum over the various dyadic ranges for
the $X_i$ and $Z_{ij}$, subject to (\ref{5.1a}) and (\ref{5.1b}).  We begin 
by considering
the summation over $Z_{ij}$.  The values of $Z_{ij}$ are powers of 2,
subject to the constraints (\ref{5.1b}).  These imply
that $F^2\ll X\ll B^4$.  Thus there are $O((\log B)^6)$ possible sets of
values for the various $Z_{ij}$.  Moreover there are $O((\log B)^5)$
sets of values for each given value of $F$.  Since $F$ runs over
powers of 2, subject to $F\ll X^{1/2}$ we conclude that
\[\sum_{Z_{ij}}F^{2/9}\ll X^{1/9}(\log B)^5.\]
We therefore deduce that
\begin{equation}\label{5.6}
\sum_{Z_{ij}}\sigma\ll (\log B)^{6}.
\end{equation}
It remains to consider the summation over values of the $X_i$, which
also run over powers of 2.  Here we
observe that
\[X^{1/4}(X_1/X_2)^{1/8}=X_1^{3/8}X_2^{1/8}X_3^{1/4}X_4^{1/4},\]
and that (\ref{5.1a}) yields
\[\sum_{X_i}X_i^e\ll B^e\]
if $e>0$.  It follows that
\[\sum_{X_1,X_2,X_3,X_4}X_1^{3/8}X_2^{1/8}X_3^{1/4}X_4^{1/4}\ll B.\]
Similarly we have
\[(X_1/X_2)^{1/4}X_3\ll X_1^{1/4}X_2^{1/4}X_3^{1/4}X_4^{1/4},\]
by (\ref{5.1}) and (\ref{5.2}), so that
\[\sum_{X_1,X_2,X_3,X_4}(X_1/X_2)^{1/4}X_3\ll B.\]
Combining this with (\ref{5.6})
completes the proof of the upper bound in our theorem.

Mathematical Institute,

24-29, St. Giles',

Oxford OX1 3LB

England

\mbox{}

email: rhb@maths.ox.ac.uk

\end{document}